\input amstex
\input epsf
\documentstyle{amsppt}
	\magnification=1200
	\rightheadtext{random walks, Liouville's theorem, and circle packings}
	\hcorrection{0.25in}

\define\cp{circle packing}
\define\lb{labelled complex}
\define\dcp{discrete complex polynomial}

\define\compk{\Bbb K}
\define\mcompk{$\compk$}
\define\comp#1{\Bbb{#1}}
\define\mcomp#1{$\comp #1$}
\define\simp#1#2{{\comp #1}^{(#2)}}
\define\msimp#1#2{$\simp #1#2$}

\define\vcompk#1{{\compk}_{#1}}
\define\mvcompk#1{$\vcompk #1$}
\define\vsimpk#1#2{\vcompk #1^{(#2)}}
\define\mvsimpk#1#2{$\vsimpk #1#2$}

\define\rplus{(0,\infty)}
\define\tini{t\in I}
\define\mtini{$\tini$}
\define\hatro{\hat \rho}
\define\mhatro{$\hatro$}

\define\inter#1{\text{\rm{int}}\, {\comp #1}^{(0)}}
\define\minter#1{$\inter #1$}
\define\bd#1{\text{\rm{bd}}\, {\comp #1}^{(0)}}
\define\mbd#1{$\bd #1$}

\define\interv#1{\text{\rm{int}}\, \vcompk #1^{(0)}}
\define\minterv#1{$\interv #1$}
\define\bdv#1{\text{\rm{bd}}\, \vcompk #1^{(0)}}
\define\mbdv#1{$\bdv #1$}

\define\p#1{\Cal{#1}}
\redefine\mp#1{$\p #1$}
\define\vp#1#2{\p #1_{#2}}
\define\mvp#1#2{$\vp #1#2$}

\define\smap#1{S_{\p #1}}
\define\msmap#1{$\smap #1$}

\define\carr#1{\text{\rm{carr}} (\p #1)}
\define\mcarr#1{$\carr #1$}


\define\br#1{\text{\rm{br}}({\p #1})}
\define\mbr#1{$\br #1$}

\def\today{\number\day\space\ifcase\month\or January \or
February \or March \or April \or May \or June \or July \or
August \or September \or October \or November \or
December\fi\space\number\year}


\topmatter
\title 
Recurrent random walks, Liouville's theorem, and circle packings
\endtitle
\author
Tomasz Dubejko
\endauthor
\affil 
Mathematical Sciences Research Institute\endaffil
\address
Mathematical Sciences Research Institute, Berkeley, CA 94720, U.S.A.\endaddress
\email
tdubejko$\@$msri.org\endemail
\thanks
Research at MSRI is supported in part by grant no.DMS-9022140.\endthanks
\abstract
It has been shown that univalent \cp s filling in the complex plane $\bold C$ are unique up to similarities of $\bold C$.
Here we prove that bounded degree branched \cp s properly covering $\bold C$ are uniquely determined, up to similarities of $\bold C$, by their branch sets.
In particular, when branch sets of the packings considered are empty we obtain the earlier result.

We also establish a \cp\ analogue of Liouville's theorem: if $f$ is a \cp\ map whose domain packing is infinite, univalent, and has recurrent tangency graph, then the ratio map associated with $f$ is either unbounded or constant.
\endabstract
\endtopmatter

\document
	\baselineskip=12pt
	\parindent=10pt
	\parskip=0pt

\heading 1~Introduction \endheading

Over the period of last several years various results linking \cp\ with other, classical, fields of mathematics have been established ([Bo],[D1],[HSc],[RS], [St1],[Th2]).
The most prominent connection is the one with analytic function theory, which was originally suggested by Thurston in [Th1].
However, recent developments ([BeSc],[HSc],[Mc]) point to an equally interesting relation between the theory of \cp\ and graph theory.
This relationship was first observed by Stephenson ([St1],[St2]) who used Markov processes and electric network-type arguments to prove Thurston's conjecture (cf. [RS],[Th1]).
The ``analytical'' side of \cp\ theory has predominantly contained, until very recently, results about univalent \cp s.
These include the existence and uniqueness statements, and the finite Riemann mapping theorem.
On the other hand, the ``graph-theoretical'' side gives, for example, an answer to the type problem for graphs in the language of \cp.

In [D1] and [D2] the author expanded ``analytic'' \cp\ theory to the non-univalent case by proving results that have well-known parallels in the classical setting.
Specifically, notions of discrete Blaschke products and discrete complex polynomials, and associated with them branched \cp s, have been introduced.
As the new concepts have arrived, new questions have been posed.
In particular, since it has been shown ([RS],[Sc],[St3]) that univalent \cp s which cover the complex plane $\bold C$ are unique up to similarities of $\bold C$, it is natural to ask whether the analogous fact holds for branched \cp s.
More precisely, is it true that \cp s which are proper branched coverings of $\bold C$ are uniquely determined, up to  similarities of $\bold C$, by their branch sets?
This paper gives the positive answer to the above question in the case of \cp s of bounded degree; the exact statement is contained in Theorem~4.2.
The proof of Theorem~4.2 uses both analytical developments of [D2] and graph-theoretical developments of [HSc] and [Mc].
With the help of [D2] we show that two proper branched \cp\ coverings of $\bold C$ with identical branch sets must have comparable radius functions; using [HSc]/[Mc] and random walk techniques we prove that there are no non-trivial bounded perturbations of a \cp\ whose combinatorial pattern of tangencies is a recurrent planar graph (see Corollary~3.2).
These two results yield the uniqueness of branched \cp s.
Particularly, we obtain the uniqueness of univalent packings.

Although this paper was motivated by studies of \cp s, it is worth mentioning that Corollary~3.2 is actually a special case of Theorem~3.1 when \lb es of the theorem are replaced by \cp s.
Roughly speaking, Theorem~3.1 says that if two labels $\rho_1$ and $\rho_2$, of a 2-complex \mcompk\ with recurrent 1-skeleton, are such that the curvature of the \lb\ $\compk(\rho_1)$ is no larger then that of $\compk(\rho_2)$ then either $\rho_1\equiv c\rho_2$ for some constant $c>0$ or $\rho_1$ is not comparable with $\rho_2$.
To prove the theorem we use variational methods (cf. [CdV]) and results about random walks.
Also, as a special case of Theorem~3.1 we obtain a \cp\ version of Liouville's theorem (Corollary~3.3).

We finish our introduction with the following remark.
As graphs can be divided into two types of classes, recurrent and transient, similarly \cp s can be divided into two types of classes, parabolic and hyperbolic (see [BSt1]); in the case of packings of bounded degree, a \cp\ is of parabolic/hyperbolic-type if and only if its combinatorial pattern of tangencies is a recurrent/transient graph, respectively.
Thus, this paper deals with \cp s of parabolic-type.
We would like to mention that the existence and uniqueness questions for hyperbolic-type branched \cp s are addressed in [D4].
Techniques used in [D4] are quite different from the ones used here as they reflect the distinction between hyperbolic-type packings and parabolic-type packings.

\heading 2~Vital Facts \endheading

In this section we will introduce essential terminology and definitions.
We will also recall several facts regarding \lb es, \cp s, and random walks.

We begin with \lb es and \cp s; the interested reader should see [Bo], [BSt2], [D1], and [D2] for more details.
Let \mcompk\ be a simplicial 2-complex given by a triangulation of a domain in $\bold C$.
We will assume, throughout this paper, that \mcompk\ is either {\it infinite without boundary} or {\it finite}, and has an orientation induced from $\bold C$.
Denote by \msimp K0, \minter K, \mbd K, \msimp K1, and \msimp K2 the sets of vertices, interior vertices, boundary vertices, edges, and faces of \mcompk, respectively.

A function $\rho: \simp K0 \to \rplus$ will be called a {\it label} for \mcompk, and $\compk(\rho)$ a {\it \lb } with label $\rho$ (cf. [Bo],[BSt2],[Th2]).
If \mcomp L is a subcomplex of \mcompk\ and $\rho$ is a label for \mcompk\ then $\comp L(\rho)$ will denote \mcomp L with the label $\rho|_{\simp L0}$.
For a face $\triangle = \langle u, v, w \rangle$ in \mcompk\ let $\alpha_{\rho}(v,\triangle)$ be given by 
$$
\alpha_{\rho}(v,\triangle):= \arccos\left[ \frac{\bigl(\rho(v)+\rho(u) \bigr)^2 + \bigl(\rho(v)+\rho(w) \bigr)^2 - \bigl(\rho(u)+\rho(w) \bigr)^2}{2\bigl(\rho(v)+\rho(u) \bigr)\bigl(\rho(v)+\rho(w) \bigr)} \right].\tag{2.1}
$$
Figure~1 shows the geometric interpretation of (2.1).

\midinsert
\epsfysize=4.0truecm
\centerline{\epsffile[156 284 460 508]{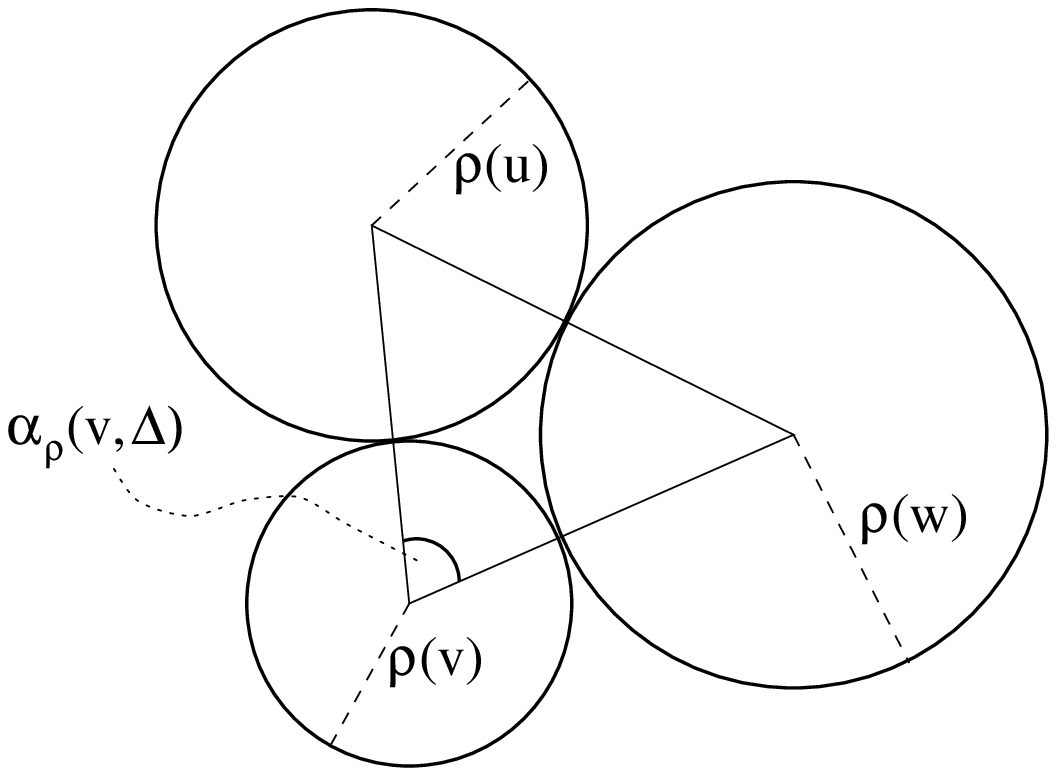}}
\captionwidth{16pc}
\botcaption{Figure~1} The angle $\alpha_{\rho}(v,\triangle)$.
\endcaption
\endinsert

If $v\in \inter K$, the quantity $\Theta_{\rho}(v):= \sum_{\triangle \in \simp K2} \alpha_{\rho}(v,\triangle)$ will be termed the {\it angle sum} of $\compk(\rho)$ at $v$, and $\Theta_{\rho} : \inter K \to \rplus$ the angle sum function.

Labelled complexes can be regarded as generalizations of more rigid structures called \cp s.
We say that a collection \mp P of circles in $\bold C$ is a {\it \cp } for \mcompk\ if there exists 1-to-1 relationship between \msimp K0 and \mp P ($\simp K0 \ni v \leftrightarrow C_{\p P}(v) \in \p P$) such that $\langle C_{\p P}(u), C_{\p P}(v), C_{\p P}(w)  \rangle$ is a positively oriented triple of mutually externally tangent circles whenever $\langle u, v, w \rangle \in \simp K2$ is a positively oriented triple.

\midinsert
\centerline{\epsfysize=7.5truecm \epsffile[82 160 547 638]{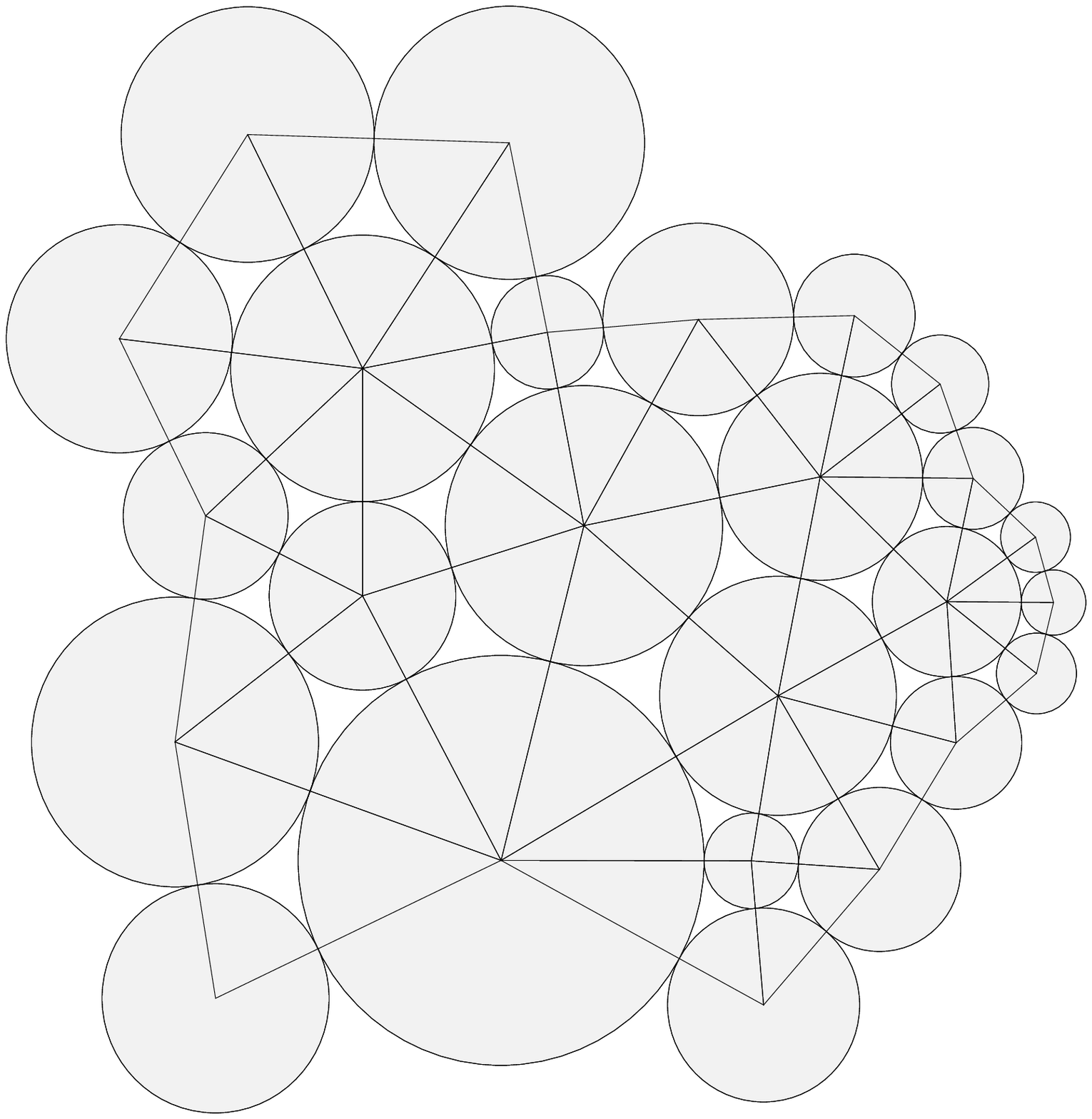} \hskip0.8truecm \epsfysize=7.5truecm \epsffile[74 150 521 637]{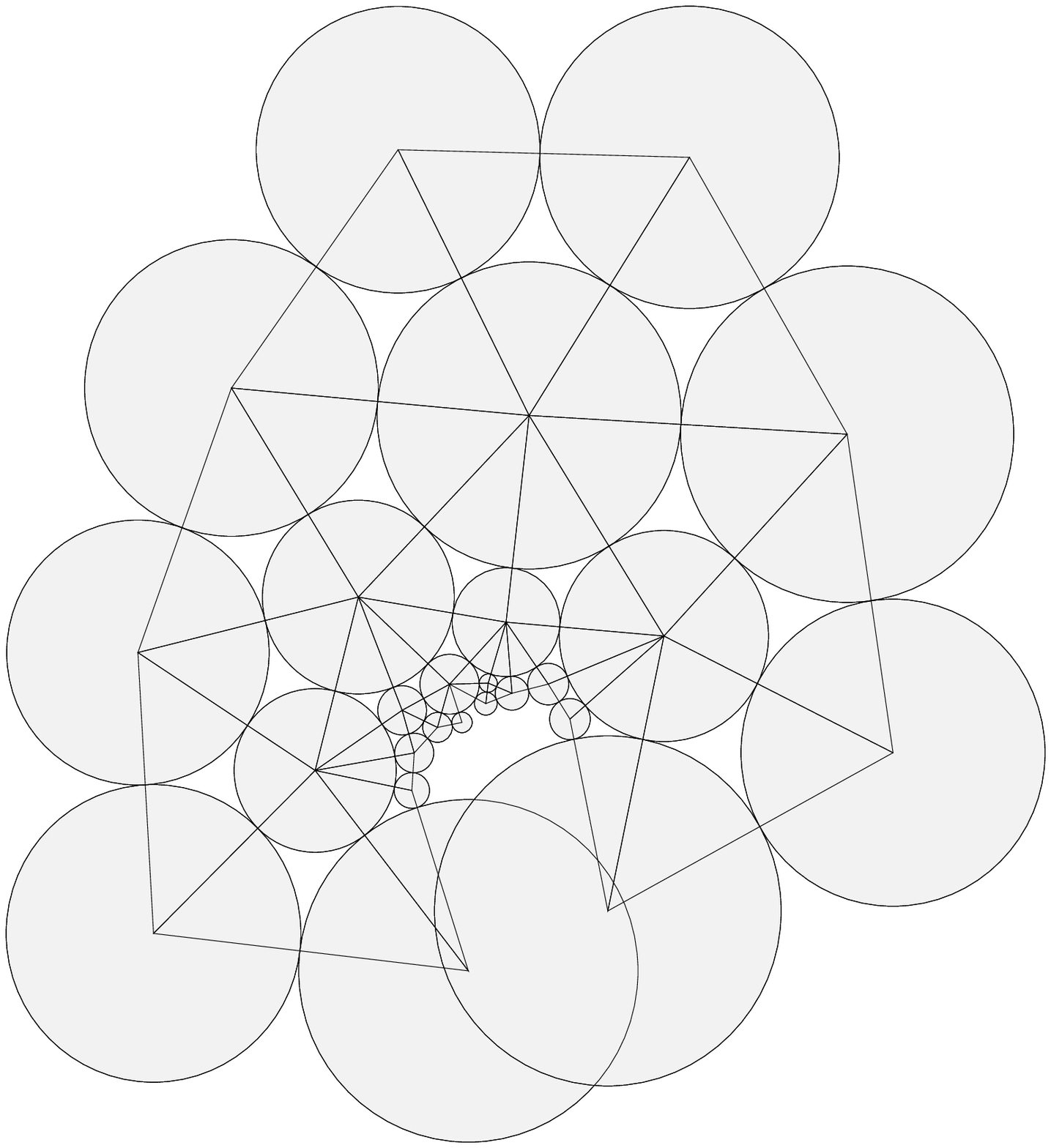}}
\captionwidth{16pc}
\botcaption{Figure~2\footnotemark"*"} Different circle packings (and their carriers) for the same 2-complex.
\endcaption
\endinsert
\footnotetext"*"{This figure was created with the help of {\tt CirclePack} by the kind permission of Ken Stephenson}%

Suppose \mp P is a \cp\ for \mcompk.
A function $r_{\p P}: \simp K0 \to \rplus$ with its value at $v\in \simp K0$ equal to the radius of $C_{\p P}(v)$ will be called the {\it radius function} of \mp P.
The {\it carrier} of \mp P, \mcarr P, is the underlying geometric complex of \mp P which is isomorphic to \mcompk\ (Figure~2).
The isomorphism is defined as the function $\smap P : \compk \to \bold C$ which maps each $v\in \simp K0$ to the center of $C_{\p P}(v)$ and then extends affinely to \msimp K1 and \msimp K2.
Hence, $\carr P := \smap P(\compk)$.

If $\triangle \in \simp K2$ and $v$ is a vertex of $\triangle$ then we define $\alpha_{\p P}(v,\triangle):= \alpha_{r_{\p P}}(v,\triangle)$.
From the earlier geometric interpretation it follows that $\alpha_{\p P}(v,\triangle)$ is the angle at the vertex $\smap P(v)$ in the euclidean triangle $\smap P(\triangle)$.
For $v\in \simp K0$, $\Theta_{\p P}(v):= \sum_{\triangle \in \simp K2} \alpha_{\p P}(v,\triangle)$ is said to be the {\it angle sum} of $\p P$ at $v$.
It follows from our definition of the \cp\ that for each $v\in \inter K$ there exists a non-negative integer $n_{\p P}(v)$ such that $\Theta_{\p P}(v)= 2(n_{\p P}(v)+1)\pi$.
If $n_{\p P}(v)>0$ then $v$ will be called a {\it branch point} of \mp P of {\it order} $n_{\p P}(v)$; equivalently we can say that the chain of circles in \mp P associated with the adjacent vertices of $v$ in \msimp K0 winds ($n_{\p P}(v)+1)$-times around $C_{\p P}(v)$.
The set $\br P\subset \inter K\times \bold N$ will denote the {\it branch set} of \mp P, i.e. $\br P = \bigl\{ (v_1, n_{\p P}(v_1)),\dots , (v_m, n_{\p P}(v_m)) \bigr\}$, where the $v_i$'s, $i=1,\dots ,m$, are all the branch points of \mp P.

It is important to realize (see [BSt2]) that if all values of the angle function $\Theta_{\rho}$ of a \lb\ $\compk(\rho)$ are positive integer multiples of $2\pi$, then $\rho\equiv r_{\p P}$ and $\Theta_{\rho}\equiv \Theta_{\p P}$ for some \cp\ \mp P for \mcompk.

When working with \lb es it is very helpful to have the notion of subpackings (see [BSt2]).
Given a function $\Theta : \inter K \to \rplus$ we say that a \lb\ $\compk(\rho)$ is a {\it subpacking} or {\it packing} for $\Theta$ if $\Theta_{\rho} \ge \Theta$ or $\Theta_{\rho} \equiv \Theta$, respectively.
In particular, a \cp\ \mp P for \mcompk\ is a subpacking or packing for $\Theta$ if the \lb\ $\compk(r_{\p P})$ is a subpacking or packing for $\Theta$, respectively.
It turns out that subpackings have a nice monotonicity property: if $\compk(\rho_1)$ and $\compk(\rho_2)$ are subpackings for $\Theta$ then $\compk(\hat \rho)$ is also a subpacking for $\Theta$, where $\hat \rho(v):=\text{max}\{\rho_1(v), \rho_2(v)\}$.
The following simple observation is known as Maximum Principle (cf. [D1]).

\proclaim{Fact~2.1{\rm (Maximum Principle)}}
Let $\compk(\rho_1)$ and $\compk(\rho_2)$ be \lb es for \mcompk.
If $\compk(\rho_2)$ is a subpacking for $\Theta_{\rho_1}$ then either $\tfrac{\rho_2(v)}{\rho_1(v)} < \sup_{v\in \simp K0}\tfrac{\rho_2(v)}{\rho_1(v)}$ for all $v\in \inter K$ or there exists $c>0$ such that $\rho_2 \equiv c\rho_1$.
In particular, if $\Theta_{\rho_1}\equiv \Theta_{\rho_2}$ and \mcompk\ is finite then $\min_{v\in \bd K}\tfrac{\rho_2(v)}{\rho_1(v)} \le \tfrac{\rho_2(v)}{\rho_1(v)} \le \max_{v\in \bd K}\tfrac{\rho_2(v)}{\rho_1(v)}$ for all $v\in \simp K0$.
\endproclaim

The next result deals with a boundary value problem for \lb es.
Its first part can be easily verified using Maximum Principle and Perron's methods introduced in [BSt2] (cf. [Bo],[D1],[G]); the second part can be found in [D1] (cf. [Bo]).

\proclaim{Theorem~2.2}
Let $\Theta : \inter K \to \rplus$ be a function.
Suppose a \lb\ $\compk(\rho)$ is a packing for $\Theta$.
Then 
\roster
\item"(i)"
For any function $\lambda : \bd K \to \rplus$ there exists a unique label $\rho_{\lambda} : \simp K0 \to \rplus$ such that $\compk(\rho_{\lambda})$ is a packing for $\Theta$ and $\rho_{\lambda}|_{\bd K} \equiv \lambda$.
\item"(ii)"
If \mcompk\ is a finite triangulation of a disc and $\Theta(\inter K)\subset \{2n\pi\}_{n>0}$ then there exists a \cp\ \mp P for \mcompk\ contained in the unit disc such that $\Theta_{\p P}\equiv \Theta$ and all boundary circles of \mp P are internally tangent to the unit circle.
Moreover, \mp P, called {\sl Bl-type packing} for \mcompk, is unique up to M\"obius transformations preserving the unit disc.
\endroster
\endproclaim

Let \mp P and \mp Q be \cp s for \mcompk.
A function $f: \carr P \to \bold C$ is said to be the {\it \cp\ map} (shortly, cp-map) from \mp P to \mp Q if $f$ maps the center of $C_{\p P}(v)$ to the center of $C_{\p Q}(v)$ and extends affinely to faces of \mcarr P.
In such a case we call \mp P (\mp Q) the {\it domain} ({\it range}) packing of $f$.
If \mp P is {\it univalent}, i.e. all circles in \mp P have mutually disjoint interiors (in particular, $\br P = \emptyset$), then \mcarr P can naturally be regarded as a subset of $\bold C$ and the map $f$ as a function defined in such a subset.
If \mp P is not univalent then $f$ should be thought of as a function from the underlying surface of \mp P (=\mcarr P).

With the pair of packings \mp P and \mp Q we associate another map $f^\# : \carr P \to \rplus$ called the {\it ratio map} from \mp P to \mp Q.
It is defined on \msimp K0 by $f^\#(\smap P(v)):=\tfrac{r_{\p Q}(v)}{r_{\p P}(v)}$ and then extended affinely to \msimp K1 and \msimp K2.
As before, depending on whether \mp P is univalent or not, $f^\#$ can be regarded as a function from a subset of $\bold C$ or from the underlying surface of \mp P.
Due to various approximation results (see [D1],[D2],[RS],[St1]), $f^\#$ should be thought of as a discrete analogue of the absolute value of the complex derivative of the cp-map $f$.

Before we state the next result we need the following three definitions.
We will say that a complex \mcompk\ is of {\it bounded degree} $d$ if every vertex in \mcompk\ has at most $d$ neighboring vertices.
A function $f:\Omega \to \bold C$, $\Omega\subset \bold C$, is of {\it valence} $M$ if $M$ is the least upper bound on the number of elements in $f^{-1}(z)$ for every point $z\in \bold C$.
Finally, the {\it valence} of a \cp\ \mp P is defined as the valence of the cp-map from a univalent \cp\ to the packing \mp P.

\proclaim{Theorem~2.3}
Let \mp P be a univalent \cp\ for \mcompk\ such that $\carr P=\bold C$ and \mcompk\ is of bounded degree $d$.
Let $M$ be a positive integer.
There exists a constant $\kappa = \kappa (M,d)$, depending only on $M$ and $d$, with the following property:
if \mp Q is a finite valence \cp\ for \mcompk\ with finite (possible empty) branch set, say $\{(v_1,n_1),\dots , (v_m,n_m)\}$, satisfying $\tsize \sum_{i=1}^{m} n_i \le M$ then:
\roster
\item
$\tfrac{r_{\p Q}(v)}{r_{\p Q}(w)} \in (1/\kappa, \kappa)$ for all pairs of adjacent vertices $v,w\in \simp K0$,
\item
the cp-map $f$ from \mp P to \mp Q is $\kappa$-quasiregular with a decomposition $f=\psi \circ h$, where $h:\bold C \to \bold C$ is $\kappa$-quasiconformal and $\psi$ is a complex polynomial with the branch set $\bigl\{ \bigl(h(S_{\p P}(v_1)),n_1\bigr),\dots , \bigl(h(S_{\p P}(v_m)),n_m\bigr)\bigr\}$.
\endroster
\endproclaim

The above fact is a result of Corollary~4.9 and Lemma~5.2 of [D2].
Notice that if $\br Q = \emptyset$ then \mp Q is univalent and $f$ is $\kappa$-quasiconformal.

We now shift our attention to the subject of random walks.
We will establish basic terminology and notation; the reader should refer to [So] and [W] for more information.
If \mcompk\ is an infinite simplicial 2-complex given by a triangulation of a domain in $\bold C$ then its 1-skeleton, \msimp K1, is a graph.
If \mcompk\ is of bounded degree then so is the graph \msimp K1.
A function $p: \simp K0\times \simp K0 \to [0,\infty)$ such that $\tsize \sum_{w} p(v,w) =1$ for each $v\in \simp K0$ will be called the {\it transition probability function} and $p(v,w)$ the {\it transition probability} from $v$ to $w$.
A stochastic process on \msimp K1 given by $p$ will be called the {\it random walk} on \msimp K1 with the transition probability $p$, and denoted $(\simp K1, p)$.
For the purposes of this paper we assume that
\roster
\item
$p(v,v)=0$ for all $v\in \simp K0$,
\item
$p(v,w)=0$ if $v$ and $w$ are not neighbors,
\item
$p(v,w)>0$ whenever $v$ and $w$ are adjacent.
\endroster

An example of a random walk on \msimp K1 is the {\it simple random walk} on \msimp K1 given by 
$$
p(v,w)=\cases
\tfrac{1}{\delta(v)}& \qquad\text{if $w$ is a neighbor of $v$}\\
	0& \qquad\text{otherwise}, \endcases
$$
where $\delta(v)$ denotes the number of adjacent vertices of $v$ in \msimp K1.

A random walk $(\simp K1, p)$ is called {\it reversible} if there exists a function $c:\simp K0 \to \rplus$ such that $c(v)p(v,w)= c(w)p(w,v)$ for all $v,w\in \simp K0$.
Notice that $c(v,w)$, called the {\it conductence} between $v$ and $w$, given by $c(v,w):=c(v)p(v,w)$ is symetric.
Conversely, a symetric function $c: \simp K0 \times \simp K0 \to \rplus$ such that $\sum_w c(v,w) \in \rplus$ for all $v\in \simp K0$ defines the reversible random walk $(\simp K1, p)$ by $p(v,w)=c(v,w)/c(v)$.
Note that in the case of the simple random walk, $c(v,w)=1$ for all neighboring vertices $v,w$.
If $f: \simp K0 \to \bold R$ is a function and $(\simp K1, p)$ is a random walk then $f$ is said to be {\it superharmonic} with respect to $p$ if for every $v\in \simp K0$
$$
\sum_{w}p(v,w)f(w) \le f(v).
$$
A random walk is called {\it recurrent} if the probability of visiting every vertex infinitely many times is equal to 1.
If a random walk is not recurrent then it is called {\it transient}.
If the simple random walk on a graph is recurrent (transient) then the graph is called recurrent (transient).
The following fact is a consequence of Theorem~4.1 and Corollary~4.14 of [W].

\proclaim{Theorem~2.4}
\roster
\item
A random walk $(\simp K1, p)$ is recurrent if and only if there are no non-constant, non-negative superharmonic functions for $p$.
\item
Let $(\simp K1, p)$ and $(\simp K1, \tilde p)$ be reversible random walks with conductances $c$ and $\tilde c$, respectively.
If $\tilde c(v,w)\le c(v,w)$ for all $v,w\in \simp K0$ then the recurrence of $(\simp K1, p)$ implies the recurrence of $(\simp K1, \tilde p)$.
\endroster
\endproclaim

We finish this section with a result that determines a connection between the type of a graph and the type of a univalent \cp\ associated with it (see [Mc], [HSc]).

\proclaim{Theorem~2.5}
Let \mcompk\ be a triangulation of an open disc.
If \msimp K1 is recurrent then there exists an infinite univalent \cp\ for \mcompk\ whose carrier is $\bold C$.
Conversely, if there exists a \cp\ for \mcompk\ whose carrier is $\bold C$ and \mcompk\ is of bounded degree then \msimp K1 is recurrent.
\endproclaim

\heading 3~Variations of labelled complexes \endheading

Our goal in this section is to show that there are no non-trivial, angle-sum-preserving, bounded perturbations of a \cp\ for an infinite 2-complex with recurrent 1-skeleton.
The precise statement is Corollary~3.2 to the following theorem.

\proclaim{Theorem~3.1}
Let \mcompk\ be a triangulation of an open domain in $\bold C$ with recurrent \msimp K1.
Let $\compk(\rho)$ and $\compk(\tilde \rho)$ be \lb es for \mcompk.
Suppose $\compk(\tilde \rho)$ is a subpacking for $\Theta_{\rho}$.
Then the function $A: \simp K0 \to \rplus$, $A(v):=\tfrac{\tilde \rho(v)}{\rho(v)}$, is bounded if and only if it is constant.
\endproclaim

A proof of the above result will be given shortly, first we draw some corollaries.

\proclaim{Corollary~3.2}
Let \mcompk\ be as in Theorem~3.1.
If \mp P and \mp Q are \cp s for \mcompk\ with $\Theta_{\p P}\equiv \Theta_{\p Q}$ then either both ratio functions from $\p P$ to $\p Q$ and from $\p Q$ to $\p P$ are unbounded or are constant.
\endproclaim

The next result is a discrete analogue of Liouville's theorem for ratio functions.

\proclaim{Corollary~3.3 {\rm (Liouville's Theorem)}}
Let \mp P be a univalent \cp\ with $\carr P = \bold C$.
Suppose $f$ is a cp-map with the domain packing \mp P, and $f^\#$ is the associated ratio function.
If the 1-skeleton of the complex of \mp P is recurrent then either $f^\#$ is unbounded or is constant.
\endproclaim

From Theorem~2.5 we get

\proclaim{Corollary~3.4}
Let \mp P be a univalent \cp\ whose complex is of bounded degree and $\carr P = \bold C$.
Then a bounded ratio function with the domain packing \mp P is constant.
\endproclaim

To prove Theorem~3.1 we will need the following definition and lemma.
We use $C^l$ to denote the class of functions with continuous $l$-th derivative.

\proclaim{Definition~3.5}
A collection $\{ \rho(t) \}_{t\in I}$, $I=[0,1]$, is called a $C^l$-family of labels for \mcompk\ if the following are satisfied:
\roster
\item"(i)"
for each \mtini, $\rho(t)$ is a label for \mcompk,
\item"(ii)"
for every $v\in \simp K0$, the map $\rho_v : I \to \rplus$ given by $\rho_v(t):=\rho(t)(v)$ is $C^l$.
\endroster
\endproclaim

\proclaim{Lemma~3.6}
Let \mcomp L be a star with $m$ boundary vertices, i.e. a triangulation of a closed disc as in Figure~3.
Let $\{ \rho(t) \}_{\tini}$ be a $C^l$-family of labels for \mcomp L.
Denote by $\Theta(t)$ the angle sum at $v_0$ given by $\rho(t)$.
If $\Theta : I \to \rplus$ is non-decreasing then for each \mtini
$$
\multline
0\le \sum_{i=1}^{m}\tsize\frac{1}{\rho_0(t) + \rho_i(t)}\left[ \frac{ \sqrt{\rho_0(t) \rho_i(t) \rho_{i-1}(t)} }{ \sqrt{\rho_0(t) + \rho_i(t) + \rho_{i-1}(t)} } + \frac{ \sqrt{\rho_0(t) \rho_i(t) \rho_{i+1}(t)} }{ \sqrt{\rho_0(t) + \rho_i(t) + \rho_{i+1}(t)} } \right] \left( \frac{\rho_i'(t)}{\rho_i(t)} - \frac{\rho_0'(t)}{\rho_0(t)} \right), \\
\qquad (\rho_{m+1}\equiv \rho_1), 
\endmultline
$$
where $\rho_i : I \to \rplus$ and $\rho_i(t)=\rho(t)(v_i)$, $i=0,1,\dots ,m$.
\endproclaim

\midinsert
\epsfysize=4.5truecm
\centerline{\epsffile[221 304 396 470]{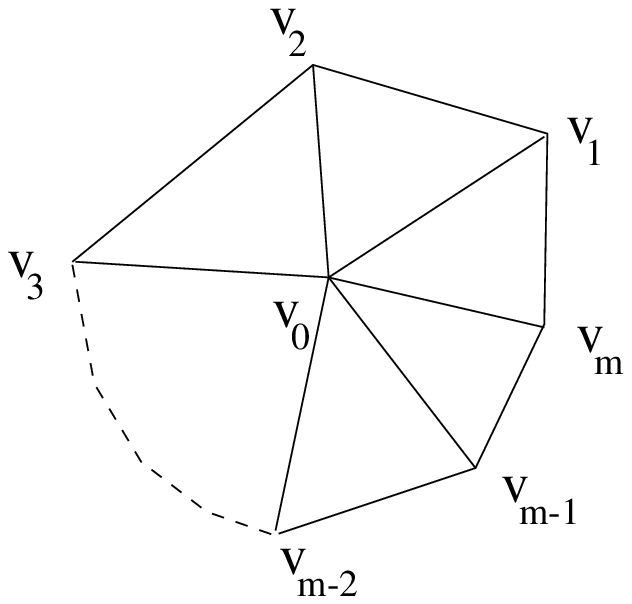}}
\captionwidth{18pc}
\botcaption{Figure~3} A star with $m$ boundary vertices.
\endcaption
\endinsert

\demo{Proof}
By the hypothesis, $0\le \Theta'(t)$ for all \mtini.
Notice that
$$
\Theta(t)=\sum_{i=1}^{m}\arccos \left[ 1-\frac{2 \rho_i(t) \rho_{i+1}(t)} {(\rho_0(t) + \rho_i(t)) (\rho_0(t) + \rho_{i+1}(t))} \right].
$$
Therefore
$$
\split
\Theta'
&=\sum_{i=1}^{m}\frac{1}{\sqrt{ 1-\left( 1-\frac{2 \rho_i \rho_{i+1}} {(\rho_0 + \rho_i) (\rho_0 + \rho_{i+1})} \right)^2 }} \left( \frac{2 \rho_i \rho_{i+1}} {(\rho_0 + \rho_i) (\rho_0 + \rho_{i+1})} \right)' 
\\
&=\sum_{i=1}^{m}\tsize\frac{1}{ \sqrt{ \frac{(\rho_0 + \rho_i) (\rho_0 + \rho_{i+1})}{ \rho_i \rho_{i+1}} -1 } }
\frac{\rho_i'\rho_0\rho_{i+1}(\rho_0 + \rho_{i+1}) + \rho_{i+1}'\rho_0\rho_i(\rho_0 + \rho_i) - \rho_0'\rho_i\rho_{i+1}(2\rho_0 + \rho_i + \rho_{i+1}) }{ \rho_i \rho_{i+1} (\rho_0 + \rho_i)(\rho_0 + \rho_{i+1})} 
\\
&=\sum_{i=1}^{m}\tsize\frac{\rho_0}{ \sqrt{ \frac{(\rho_0 + \rho_i)(\rho_0 + \rho_{i+1})}{ \rho_i \rho_{i+1}} -1 } }
\left[ \frac{\rho_i'}{\rho_i} \frac{1}{\rho_0 + \rho_i} + \frac{\rho_{i+1}'}{\rho_{i+1}} \frac{1}{\rho_0 + \rho_{i+1}} - \frac{\rho_0'}{\rho_0}\left( \frac{1}{\rho_0 + \rho_i} + \frac{1}{\rho_0 + \rho_{i+1}} \right) \right] 
\\
&=\sum_{i=1}^{m}\left[ \frac{1}{ \sqrt{ \frac{(\rho_0 + \rho_i)(\rho_0 + \rho_{i-1})}{ \rho_i \rho_{i-1}} -1 } } + \frac{1}{ \sqrt{ \frac{(\rho_0 + \rho_i)(\rho_0 + \rho_{i+1})}{ \rho_i \rho_{i+1}} -1 } } \right] \frac{\rho_0}{\rho_0 + \rho_i} \left( \frac{\rho_i'}{\rho_i} - \frac{\rho_0'}{\rho_0} \right) 
\\
&=\sum_{i=1}^{m}\frac{1}{\rho_0 + \rho_i}\left[ \frac{ \sqrt{\rho_0 \rho_i \rho_{i-1}} }{ \sqrt{\rho_0 + \rho_i + \rho_{i-1}} } + \frac{ \sqrt{\rho_0 \rho_i \rho_{i+1}} }{ \sqrt{\rho_0 + \rho_i + \rho_{i+1}} } \right] \left( \frac{\rho_i'}{\rho_i} - \frac{\rho_0'}{\rho_0} \right) \ge 0
\\
& \hphantom{=\sum_{i=1}^{m}\frac{-1}{\rho_0 + \rho_i}\left[ \frac{ \sqrt{\rho_0 \rho_i \rho_{i-1}} }{ \sqrt{\rho_0 + \rho_i + \rho_{i-1}} } + \frac{ \sqrt{\rho_0 \rho_i \rho_{i+1}} }{ \sqrt{\rho_0 + \rho_i + \rho_{i+1}} } \right] \left( \frac{\rho_i'}{\rho_i} - \frac{\rho_0'}{\rho_0} \right) \ge 0.
}
\quad\qed
\endsplit
$$
\enddemo

We can now establish Theorem~3.1.
The proof will proceed by contradiction; using variational methods we will generate a random walk which, by the hypothesis of the theorem, cannot exist.
The reader might be interested to learn that Colin de Verdi\`ere have also applied variational techniques in [CdV] to show the existence of \cp s.

\demo{Proof of Theorem~3.1}
Suppose that the assertion is false.
Then there exist two vertices $w_0, w_1\in \simp K0$ such that $\tfrac{\tilde \rho (w_0)}{\rho (w_0)}\ne \tfrac{\tilde \rho (w_1)}{\rho (w_1)}$ and $M:=\sup_{v\in \simp K0} A(v) < \infty$.
We may assume without loss of generality that $\tilde \rho (w_0) = \rho (w_0)$ and $\tilde \rho (w_1) > \rho (w_1)$.
Define a family $\bold F$ of labels for \mcompk\ by
$$
\bigl\{ \text{$\varrho$: $\compk(\varrho)$ is a subpacking for $\Theta{\rho}$, $\varrho (w_0) = \rho (w_0)$, and $\tfrac{\varrho (v)}{\rho (v)}\le M$ for all $v\in \simp K0$}\bigr\}.
$$
Notice that $\bold F$ is not empty as $\tilde \rho, \rho \in \bold F$.
It follows from a monotonicity property of subpackings that the label \mhatro\ given by $\hatro(v) := \sup_{\varrho \in \bold F}\{\varrho (v)\}$ is a subpacking for $\Theta{\rho}$.
From definitions of $\bold F$ and \mhatro\ we obtain that $\hatro\in \bold F$, and $\hatro\ge \rho$ with $\hatro (w_1)> \rho (w_1)$.
Hence \mhatro\ is not a constant multiple of $\rho$.
Moreover, $\sup_{v\in \simp K0}\tfrac{\hatro (v)}{\rho (v)} = M$.
From Maximum Principle, $\tfrac{\hatro (v)}{\rho (v)} < M$ for all $v\in \simp K0$.
This implies that $\Theta_{\hatro}(v)= \Theta_{\rho}(v)$ for all $v\ne w_0$, and $\Theta_{\hatro}(w_0) > \Theta_{\rho}(w_0)$ (for if $\Theta_{\hatro}(v_0) > \Theta_{\rho}(v_0)$ for some $v_0\ne w_0$ then there is $\epsilon >0$ such that $\bar \rho: \simp K0 \to \rplus$, defined by $\bar \rho|_{\simp K0 \setminus \{v_0\}} \equiv  \hatro|_{\simp K0 \setminus \{v_0\}}$ and $\bar \rho(v_0)=\hatro(v_0) + \epsilon$, belongs to $\bold F$).

Let $\{\vcompk n\}$ be a sequence of finite, connected 2-subcomplexes of \mcompk\ such that $\vcompk n \subset \vcompk {n+1}$, $\bigcup_{n=1}^{\infty}\vcompk n = \compk$ , and $w_0, w_1\in \interv 1$.
It follows from Theorem~2.2(i) (and Implicite Function Theorem) that for each $n$ there is a unique $C^1$-family of labels $\{\rho_n(t)\}_{\tini}$ for \mvcompk n such that 
\roster
\item
$\rho_n(t)(v)= (1-t)\rho(v) + t\hatro(v)$ for every $v\in \bdv n$,
\item
$\rho_n(t)(w_0) = \rho(w_0)$ for all $t$,
\item
$\Theta_{\rho_n(t)}(v)= \Theta_{\rho}(v)$ for every $v\in \interv n \setminus \{w_0\}$ and all $t$.
\endroster
Observe that the uniqueness in Theorem~2.2(i) gives $\rho_n(0)(v)= \rho(v)$ and $\rho_n(1)(v)= \hatro(v)$ for all $v\in \vsimpk n0$.
Since for each $v\in \bdv n$ the function $\rho_{v,n}: I \to \rplus$, $\rho_{v,n}(t):= \rho_n(t)(v)$, is non-decreasing, Maximum Principle implies that for every $v\in \vsimpk n0$ the function $\rho_{v,n}$ is non-decreasing.
Thus the function $\Theta_{w_0,n} : I \to \rplus$, $\Theta_{v,n}(t):= \Theta_{\rho_n(t)}(v)$, is non-decreasing.
Moreover, from Maximum Principle we get for all \mtini
$$
\multline
1=\min_{v\in \{w_0\}\cup \bdv n} \frac{\rho_{v,n}(t)}{\rho(v)} = 
\min_{v\in \vsimpk n0} \frac{\rho_{v,n}(t)}{\rho(v)} \le \\
\max_{v\in \vsimpk n0} \frac{\rho_{v,n}(t)}{\rho(v)} =
\max_{v\in \bdv n} \frac{\rho_{v,n}(t)}{\rho(v)} \le
\max_{v\in \bdv n} \frac{\hatro(v)}{\rho(v)} < M.
\endmultline
$$
Hence 
$$
\rho_{v,n}(t) \in [\rho(v),M \rho(v)] \qquad \text{for every $n$ and all \mtini}.
\tag"(3.1)"
$$
Notice that using Maximum Principle we also obtain the following inequalities for small $h>0$
$$
\multline
0\le \rho_{v,n}(t+h)- \rho_{v,n}(t) = 
\rho_{v,n}(t) \left( \tfrac{\rho_{v,n}(t+h)}{\rho_{v,n}(t)} -1 \right) \le \\
\rho_{v,n}(t) \max_{v\in \bdv n} \left( \tfrac{\rho_{v,n}(t+h)}{\rho_{v,n}(t)} -1 \right) =
\rho_{v,n}(t) \max_{v\in \bdv n} \left( \tfrac{\rho_{v,n}(t+h) - \rho_{v,n}(t)}{\rho_{v,n}(t)} \right).
\endmultline
$$

Thus, by dividing both sides in the above inequalities by $h$ and then taking the limit as $h\to 0$ we get
$$
0\le \frac{\rho_{v,n}'(t)}{\rho_{v,n}(t)}\le 
\max_{v\in \bdv n} \frac{\rho_{v,n}'(t)}{\rho_{v,n}(t)} \qquad \text{for all $v\in \vsimpk n0$, \mtini, and every $n$}.
$$
In particular,
$$
0\le \frac{\rho_{v,n}'(t)}{\rho_{v,n}(t)}\le M \qquad \text{for all $v\in \vsimpk n0$, \mtini, and every $n$}. \tag"(3.2)"
$$

Denote by $\mu:= \log \tfrac{\hatro(w_1)}{\rho(w_1)}$.
Since $\mu>0$ and 
$$
\mu= \int_0^1 \left(\log \rho_{w_1,n}(t)\right)'dt = \int_0^1 \frac{\rho_{w_1,n}'(t)}{\rho_{w_1,n}(t)} dt ,
$$
for each $n$ there exists $t_n\in I$ such that 
$$
0<\mu\le \frac{\rho_{w_1,n}'(t_n)}{\rho_{w_1,n}(t_n)}.
\tag"(3.3)"
$$

Suppose $v$ is an interior vertex in \mvsimpk n0 and $w\in \vsimpk n0$ is a neighbor of $v$.
Let $w'$ and $w''$ be the other two vertices in \mvsimpk n0 adjacent to both $v$ and $w$.
We define 
$$
c_n(v,w;t):= \tsize\frac{1}{\rho_{v,n}(t)+\rho_{w,n}(t)}
\left[ \sqrt{ \frac{\rho_{v,n}(t) \rho_{w,n}(t) \rho_{w',n}(t)}{\rho_{v,n}(t)+\rho_{w,n}(t)+\rho_{w',n}(t)} } + 
\sqrt{ \frac{\rho_{v,n}(t) \rho_{w,n}(t) \rho_{w'',n}(t)}{\rho_{v,n}(t)+\rho_{w,n}(t)+\rho_{w'',n}(t)} } \right].
\tag"(3.4)"
$$
Notice that if $w$ is also an interior vertex in \mvsimpk n0 then 
$$
c_n(v,w;t)=c_n(w,v;t).
\tag"(3.5)"
$$
Moreover
$$
c_n(v,w;t)\le 2\frac{\sqrt{\rho_{v,n}(t) \rho_{w,n}(t)}} {\rho_{v,n}(t)+\rho_{w,n}(t)} \le 1.
\tag"(3.6)"
$$
From Lemma~3.6 we get for every $n$ and $t$
$$
\aligned
0&\le \sum_w c_n(w_0,w;t)\left( \frac{\rho_{w,n}'(t)}{\rho_{w,n}(t)} - 
\frac{\rho_{w_0,n}'(t)}{\rho_{w_0,n}(t)} \right) \qquad \text{and}\\
0&= \sum_w c_n(v,w;t)\left( \frac{\rho_{w,n}'(t)}{\rho_{w,n}(t)} - 
\frac{\rho_{v,n}'(t)}{\rho_{v,n}(t)} \right) \qquad \text{for $v\in \interv n \setminus \{w_0\}$}.
\endaligned
\tag"(3.7)"
$$

We now make the following observation.
Suppose $\{s_n\}$ is a sequence of functions $s_n:\simp K0 \times I \to \bold R$ such that for each $v\in \simp K0$, the $s_n(v,\cdot): I \to \bold R$, $n=1,2, \dots$, have uniformly bounded ranges.
Then for any sequence $\{a_n\}$ of points in $I$ there are subsequences $\{s_{n(k)}\}$ and $\{a_{n(k)}\}$ such that $\lim_{k\to \infty}s_{n(k)}\bigl(v, a_{n(k)}\bigr)$ exists for every $v\in \simp K0$.
This is due to the following diagonalization argument:
\block
{\eightpoint
We begin by enumerating all vertices in \msimp K0 so that $\{v_i\}_{i=1}^{\infty} = \simp K0$.
Since the maps $s_n(v_1,\cdot)$ have uniformly bounded ranges there exist subsequences $\{s_{n_1(k)}\}$ and $\{a_{n_1(k)}\}$ of $\{s_n\}$ and $\{a_n\}$, respectively, such that $\{s_{n_1(k)}\bigl(v_1, a_{n_1(k)}\bigr)\}_{k=1}^{\infty}$ is convergent.
Next, we repeat the above construction for $v_2$ with $\{s_{n_1(k)}\}$ and $\{a_{n_1(k)}\}$ instead of $\{s_n\}$ and $\{a_n\}$, respectively.
This yields two subsequences $\{s_{n_2(k)}\}$ and $\{a_{n_2(k)}\}$ of $\{s_{n_1(k)}\}$ and $\{a_{n_1(k)}\}$, respectively, such that $\{s_{n_2(k)}\bigl(v_2, a_{n_2(k)}\bigr)\}_{k=1}^{\infty}$ is convergent and, of course, $\lim_{k\to \infty}s_{n_2(k)}\bigl(v_1, a_{n_2(k)}\bigr) =\lim_{k\to \infty}s_{n_1(k)}\bigl(v_1, a_{n_1(k)}\bigr)$.
We continue this process for all $v_i$'s.
Then by setting $n(k):=n_k(k)$ we obtain desired sequences $\{s_{n(k)}\}$ and $\{a_{n(k)}\}$.
}\endblock
We finish our observation by noting that analogous diagonalization argument is valid when \msimp K0 in the domains of the functions $s_n$ is replaced by $\simp K0\times \simp K0$.

We now return to our sequence $\{t_n\}$ satisfying (3.3).
We can assume, by using the inequalities (3.2) and (3.6) and by applying diagonalization argument if necessary, that
$$
f(v):=\lim_{n\to \infty}\frac{\rho_{v,n}'(t_n)}{\rho_{v,n}(t_n)} \quad\text{and}\quad
c(v,w):=\lim_{n\to \infty}c_n(v,w;t_n) \quad\text{exist for all $v,w\in \simp K0$}.
$$
From (3.1), (3.4), (3.5), and (3.6) it follows that
$$
0< c(v,w)=c(w,v)\le 1 \quad\text{for all adjacent vertices $v$ and $w$ in \msimp K0}.
\tag"(3.8)"
$$
Also, (3.2) and (3.3) imply
$$
0\le f(v)\le M \quad\text{for all $v\in \simp K0$, and} \quad 0=f(w_0)< \mu \le f(w_1).
\tag"(3.9)"
$$
Moreover, from (3.7) we get
$$
0\le \sum_w p(v,w)\left(f(w)-f(v)\right) \quad\text{for every $v\in \simp K0$},
$$
where $p(v,w)=c(v,w)/\sum_w c(v,w)$.
Thus the function $\tilde f: \simp K0 \to \bold R$, $\tilde f(v):= M-f(v)$, is superharmonic for the reversible random walk $(\simp K1, p)$ with conductance between $v$ and $w$ equal to $c(v,w)$.
From Theorem~2.4(1) and (3.9) we conclude that $(\simp K1, p)$ is transient.
On the other hand, since by the hypothesis \msimp K1 is recurrent, Theorem~2.4(2) and (3.8) (recall that the conductances for a simple random walk are all equal to 1) imply that $(\simp K1, p)$ is recurrent.
This is a contradiction.
Hence $\frac{\tilde \rho(w)}{\rho(w)} = \frac{\tilde \rho(v)}{\rho(v)}$ for all $v,w\in \simp K0$, and the proof is complete.
\qed
\enddemo

\heading 4~Uniqueness of branched packings \endheading

We begin this section with the following definition which was originally introduced in [D2].

\proclaim{Definition~4.1}
Let \mcompk\ be a triangulation of an open disc.
A cp-map $f$ is said to be a discrete complex polynomial for \mcompk\ if the carrier of the domain packing of $f$ is equal to $\bold C$ and the range packing is a proper, finite, branched covering of $\bold C$, i.e. $f$ has a decomposition $f=\psi \circ h$, where $h$ is a self-homeomorphism of $\bold C$ and $\psi$ is a complex polynomial.
\endproclaim

It was shown in [D2] that \dcp s exist.
The range packings of such maps form a natural expansion of the class of univalent \cp s that cover $\bold C$.
However, the uniqueness result for branched packings similar to the one for univalent packings (see [RS],[Sc],[St3]) was proved so far only in the case of packings for regular hexagonal lattice ([D2]).
Here we will show the uniqueness of branched circle packings of bounded valence for bounded degree triangulations whose 1-skeletons are recurrent; the case when 1-skeletons are transient is treated in [D4].
In particular, we obtain that range packings of \dcp s are unique, and furthermore, when branch sets of such packings are empty we get the earlier results of [RS], [Sc], and [St3], which are crucial for approximation schemes (cf. [D1],[D2],[H],[RS],[St1]).

\proclaim{Theorem~4.2}
Let \mcompk\ be a bounded degree triangulation of an open disc with recurrent \msimp K1.
Let $\bold F_{\goth B}$ denote the set of all bounded valence \cp s for \mcompk\ with identical finite branch sets equal to $\goth B$.
If $\bold F_{\goth B}\ne \emptyset$ then there exists $\p P\in \bold F_{\goth B}$ such that \mp P is the range packing of a \dcp\ and $\bold F_{\goth B}= \{z\p P\}_{z\in \bold C}$, i.e. all packings in $\bold F_{\goth B}$ are copies of \mp P under similarities of $\bold C$.
\endproclaim

\demo{Proof}
From Theorem~2.5 we have that there exists a univalent \cp\ \mp O for \mcompk\ such that $\carr O =\bold C$.
Since $\bold F_{\goth B}\ne \emptyset$ there is $\p P\in \bold F_{\goth B}$ of finite valence with branch set $\goth B$.
From Theorem~2.3 it follows that the cp-map from \mp O to \mp P is a \dcp.
As \mp P was arbitrary we get that all packings in $\bold F_{\goth B}$ are range packings of \dcp s.
We also see that $\{z\p P\}_{z\in \bold C}\subset \bold F_{\goth B}$.
To finish the proof it is sufficient to show that for each $\p Q \in \bold F_{\goth B}$ there exists a constant $c_{\p Q}$ such that $r_{\p Q}(v)= c_{\p Q} r_{\p P}(v)$ for all $v\in \simp K0$.
However, by Corollary~3.2 it is enough to show that for each $\p Q \in \bold F_{\goth B}$ there is $M_{\p Q}$ such that
$$
\frac1{M_{\p Q}}\le \frac{r_{\p Q}(v)}{r_{\p P}(v)} \le M_{\p Q}
\qquad\text{for all $v\in \simp K0$}.
\tag"(*)"
$$
The remainder of our proof is devoted entirely to verification of (*).

Fix $\p Q \in \bold F_{\goth B}$, $\p Q\ne \p P$.
By applying similarities to \mp O, \mp P, and \mp Q if necessary, we can assume that $\smap O(v_0)=\smap P(v_0)=\smap Q(v_0)=0$ and $r_{\p O}(v_0)=r_{\p P}(v_0)=r_{\p Q}(v_0)=1$ for a designated vertex $v_0\in \simp K0$, where the notation is as in Section~2.
Denote by $f_1$ and $f_2$ the cp-maps from \mp O to \mp P and from \mp O to \mp Q, respectively.
Notice that $f_1(0)=0=f_2(0)$ and, by Theorem~2.3, $f_1=\psi_1\circ h_1$ and $f_2=\psi_2\circ h_2$, where $h_1$, $h_2$ are entire $K$-quasiconformal mappings, $\psi_1$, $\psi_2$ are complex polynomials of degree, say, $l$, and $h_1(0)=0=h_2(0)$, $h_1(\infty)=\infty=h_2(\infty)$.
Moreover, for each $j$, $j=1,2$,
$$
\frac1K\le \frac{r_j(w)}{r_j(v)}\le K \quad\text{for all pairs of adjacent vertices $v,w \in \simp K0$},
\tag"(4.1)"
$$
where $r_1$ and $r_2$ are the radius functions of \mp P and \mp Q, respectively.
It follows from quasiconformal arguments (Theorem~2.4 of [L])  that there exists a constant $c=c(K)$ depending only on $K$ such that
$$
\frac{\max_{\alpha}|h_j(\delta e^{i\alpha})|}{\min_{\alpha}|h_j(\delta e^{i\alpha})|} \le c \qquad\text{for every $\delta >0$, $j=1,2$}.
$$
Since $\psi_j$'s are polynomials of degree $l$, there is $R>0$ such that for each $j$, $j=1,2$,
$$
\frac{ \max_{\alpha}|f_j(\delta e^{i\alpha})| }{ \min_{\alpha}|f_j(\delta e^{i\alpha})| } =
\frac{ \max_{\alpha}|h_j(\delta e^{i\alpha})| }{ \min_{\alpha}|h_j(\delta e^{i\alpha})| }
\frac{|1+ O(\frac1{\delta})|}{|1- O(\frac1{\delta})|}\le 2c
\quad\text{for all $\delta >R$}.
\tag"(4.2)"
$$

For each $n>R$ let \mvp On denote the largest, simply connected portion of \mp O contained in $D(n)$, $D(\delta):=\{|z|<\delta \}$.
Write $\vcompk n \subset \compk$ for the complex of \mvp On.
Let \mvp Pn and \mvp Qn be portions of \mp P and \mp Q, respectively, associated with the complex \mvcompk n.
By Theorem~2.2(ii) there exists a Bl-type packing \mvp Bn for \mvcompk n with branch set $\goth B$, and normalized so that the circle $C_{\vp Bn}(v_0)$ is centered at 0.
Define $\rho_j(n):= \max_{\alpha}|f_j(n e^{i\alpha})|$ for $j=1,2$.
Write $\tilde \vp Pn := \tfrac1{\rho_j(n)} \vp Pn$ and $\tilde \vp Qn := \tfrac1{\rho_j(n)} \vp Qn$.
Since $\tilde \vp Pn, \tilde \vp Qn \subseteq D(1)$ and packings \mvp Bn, \mvp Pn, and \mvp Qn have identical branch sets, by discrete Schwarz Lemma of [BSt2] one has
$$
\frac1{\rho_j(n)}= \frac{r_j(v_0)}{\rho_j(n)}\le \eta_n(v_0) \qquad\text{for all $n$, and $j=1,2$},
\tag"(4.3)"
$$
where $\eta_n$ is the radius function of \mvp Bn.

From (4.2) we obtain that $\tfrac1{\rho_j(n)}f_j\bigl( \partial D(n) \bigr) \subseteq \Omega:=\{ |z|\ge \tfrac1{2c}\}$.
However, it is not generally true that all boundary circles of $\tilde \vp Pn$ and $\tilde \vp Qn$ (i.e., circles in $\tilde \vp Pn$ and $\tilde \vp Qn$ corresponding to boundary vertices of \mvcompk n) intersect $\Omega$.
Nevertheless, we will now show that they are not far away from the set $\Omega$.
To be more precise, for $v\in \bdv n$ we define
$$
\delta_j^n(v)=\cases
0& \text{if $\frac1{\rho_j(n)}\left( |f_j(v)| + r_j(v) \right)\ge \frac1{2c}$ }\\
\frac1{2c} - \frac1{\rho_j(n)}\left( |f_j(v)| + r_j(v) \right)& \text{otherwise}, \endcases
$$
which is the distance between the boundary circle $C_{\tilde \vp Pn}(v)$ and $\Omega$ when $j=1$ or between the boundary circle $C_{\tilde \vp Qn}(v)$ and $\Omega$ when $j=2$ (Figure~4).
We are only interested in vertices $v\in \bdv n$ with $\delta_j^n(v)>0$; for if $\delta_j^n(v)=0$ then $C_{\tilde \vp Pn}(v)\cap \Omega \ne\emptyset$ when $j=1$ and $C_{\tilde \vp Qn}(v)\cap \Omega \ne\emptyset$ when $j=2$.

\midinsert
\epsfysize=7.5truecm
\centerline{\epsffile[130 222 478 570]{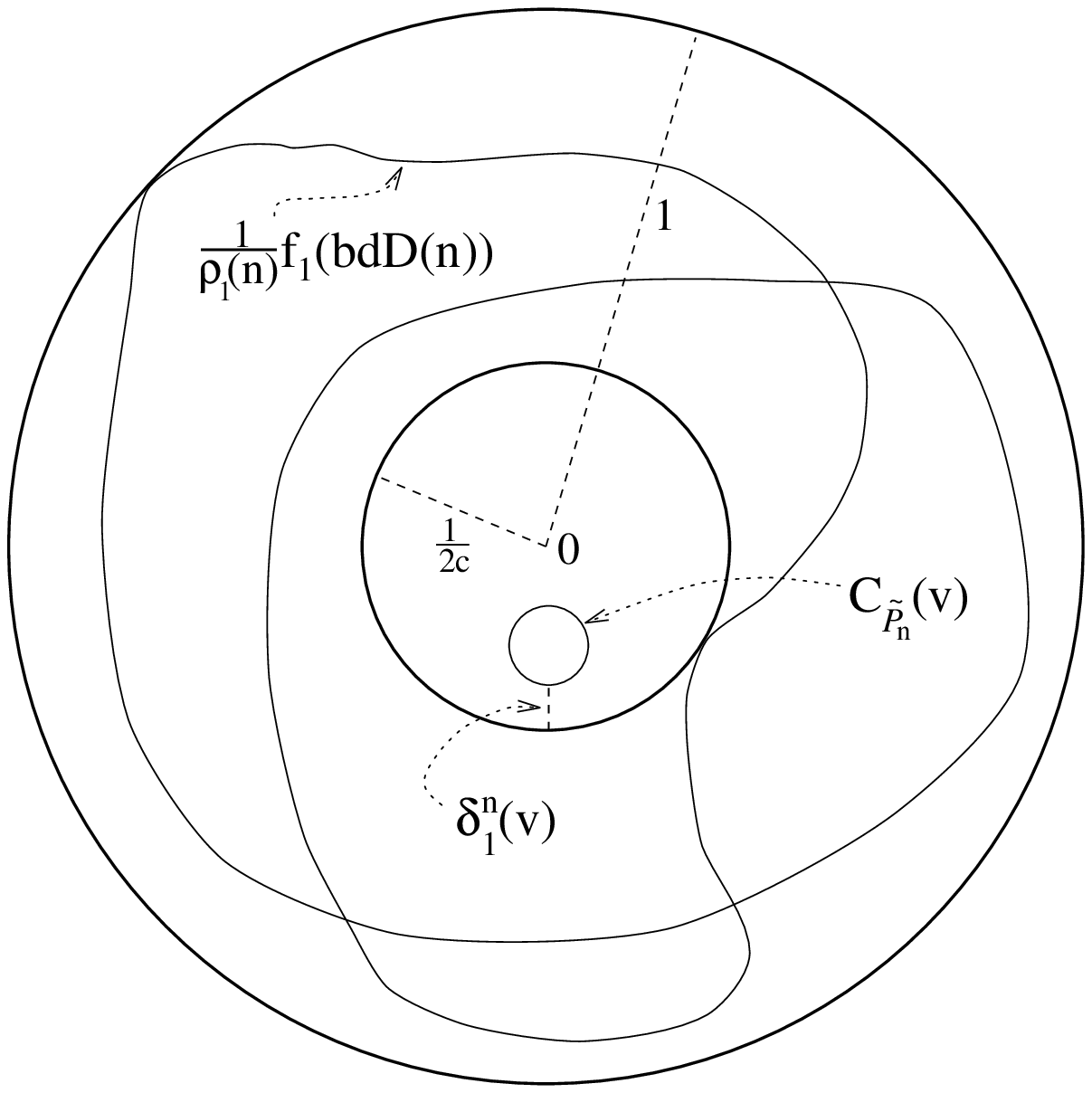}}
\captionwidth{16pc}
\botcaption{Figure~4} $\delta_j^n(v)$ for $j=1$.
\endcaption
\endinsert

It is immediate that
$$
\frac1{\rho_j(n)}r_j(v)\le \frac1{2c} - \delta_j^n(v)
\qquad\text{for $\delta_j^n(v)>0$}.
\tag"(4.4)"
$$
From the definitions of \mvcompk n and the cp-maps $f_j$ it is easy to see that if $v\in \bdv n$ and $C_{\tilde \vp Pn}(v)\cap \Omega = \emptyset$ (respectively, $C_{\tilde \vp Qn}(v)\cap \Omega = \emptyset$) then there is at least one neighbor $w\in \simp K0$ of $v$ such that the circle in $\tfrac1{\rho_1(n)} \p P$ (respectively, $\tfrac1{\rho_2(n)} \p Q$) corresponding to $w$ intersects $\Omega$.
Therefore, it follows that 
$$
\delta_j^n(v)\le \frac2{\rho_j(n)} \max\Sb w\in \simp K0\\ \text{$w$ adjucent to $v$} \endSb r_j(w) \qquad\text{whenever $\delta_j^n(v)>0$, $j=1,2$}.
$$
Thus, by (4.1) 
$$
\delta_j^n(v)\le \frac2{\rho_j(n)} Kr_j(v) \qquad\text{when $\delta_j^n(v)>0$}.
\tag"(4.5)"
$$
Hence, (4.4) and (4.5) imply 
$$
\delta_j^n(v)\le \frac{K}{c}\left( 2K+1 \right)^{-1},
\tag"(4.6)"
$$
i.e. the boundary circles of $\tilde \vp Pn$ and $\tilde \vp Qn$ that do not intersect $\Omega$ are ``not far away'' from $\Omega$.

It follows from (4.6) that $2c(2K+1)\ge (\tfrac1{2c}-\delta_j^n(v))^{-1}$.
The last inequality implies that packings $\hat \vp Pn:= 2c(2K+1)\tilde \vp Pn$ and $\hat \vp Qn:= 2c(2K+1)\tilde \vp Qn$ have their boundary circles intersecting $\{|z|\ge 1\}$, and therefore, in the language of [DSt], $\hat \vp Pn$ and $\hat \vp Qn$ properly cover $D(1)$.
Hence, by Discrete Distortion Lemma of [DSt]
$$
\eta_n(v_0)\le 2c(2K+1)\frac{r_j(v_0)}{\rho_j(n)}
\qquad\text{for all $n>R$, and $j=1,2$}.
\tag"(4.7)"
$$
Thus, (4.3) and (4.7) give
$$
\frac1{2c(2K+1)}\frac{r_1(v_0)}{\rho_1(n)}\le \frac{r_2(v_0)}{\rho_2(n)}\le 2c(2K+1) \frac{r_1(v_0)}{\rho_1(n)} \qquad\text{for all $n>R$}.
\tag"(4.8)"
$$

We will now show that (4.8) can be extended to every vertex of \mcompk\ with, perhaps, slightly larger but uniform constant instead of $2c(2K+1)$.
Denote by $\tilde \vp On$ a univalent Bl-type packing (i.e., with empty branch set) for \mvcompk n with the circle $C_{\tilde \vp On}(v_0)$ centered at 0.
Since \mcompk\ is recurrent, it follows from the proof of Discrete Uniformization Theorem of [BSt1,\S 4] and Theorem~2.5 that for each $v\in \simp K0$, $C_{\tilde \vp On}(v)$ converges to a point which is the origin as $n\to \infty$.
The monotonicity results of [BSt2, Theorem~4] imply that for each $v\in \simp K0$, $C_{\vp Bn}(v)$ converges to a point which is the origin as $n\to \infty$.

Let $\epsilon>0$ be small.
Fix $w\in \simp K0$, $w\ne v_0$.
From the above convergence results there exists $R_w\ge R$ such that 
$$
C_{\vp Bn}(w)\subset D(\epsilon), \quad
\tfrac1{\rho_1(n)}|f_j(w)|< \epsilon, \quad\text{and $\tfrac1{\rho_1(n)}|r_j(w)|< \epsilon$} \quad\text{for all $n\ge R_{w}$}.
$$
Let $\vp Bn(w)$ be a Bl-type packing for \mvcompk n with branch set $\goth B$ and normalized so that its circle associated with $w$ is centered at 0.
We also define $\tilde \vp Pn(w)$ and $\tilde \vp Qn(w)$ as packings that are obtained by translating packings $(1-\epsilon) \tilde \vp Pn$ and $(1-\epsilon) \tilde \vp Qn$, respectively, so that their circles corresponding to the vertex $w$ are centered at 0.
Thus, by discrete Schwarz Lemma of [BSt2],
$$
(1-\epsilon)\frac{r_j(w)}{\rho_j(n)}\le \eta_n^w(w) \qquad\text{for all $n\ge R_w$, and $j=1,2$},
\tag"(4.9)"
$$
where $\eta_n^w$ is the radius function of $\vp Bn(w)$.
It is clear from the definition of $\tilde \vp Pn(w)$ that the boundary circle of $\tilde \vp Pn(w)$ associated with a vertex $v\in \bdv n$ is within the distance $\delta_1^n(v) + \epsilon(\tfrac1{2c} - \delta_1^n(v)) + \epsilon$ from $\Omega$ for all $n\ge R_w$.
Similarly for boundary circles of $\tilde \vp Qn(w)$.
Thus it follows from (4.6) that $\hat \vp Pn(w):= c(K,\epsilon) \tilde \vp Pn(w)$ and $\hat \vp Qn(w):= c(K,\epsilon) \tilde \vp Qn(w)$ properly cover $D(1)$, where $c(K,\epsilon):= \tfrac{2c(2K+1)}{1- \epsilon - \epsilon 2c(2K+1)}$.
Hence, by Discrete Distortion Lemma of [DSt],
$$
\eta_n^w(w)\le c(K,\epsilon)(1-\epsilon)\frac{r_j(w)}{\rho_j(n)} \qquad\text{for all $n\ge R_w$, and $j=1,2$}.
\tag"(4.10)"
$$
Therefore, (4.9) and (4.10) give
$$
\frac1{c(K,\epsilon)} \frac{\rho_2(n)}{\rho_1(n)} \le
\frac{r_2(w)}{r_1(w)} \le c(K,\epsilon) \frac{\rho_2(n)}{\rho_1(n)} 
\qquad\text{for all $n\ge R_w$}.
\tag"(4.11)"
$$
Now (4.8) and (4.11) imply 
$$
\frac1{2c(2K+1)c(K,\epsilon)} \frac{r_2(v_0)}{r_1(v_0)} \le
\frac{r_2(w)}{r_1(w)} \le 2c(2K+1)c(K,\epsilon) \frac{r_2(v_0)}{r_1(v_0)}.
$$
Since $w$ is arbitrary vertex of \mcompk, the last double inequality proves the assertion (*).
\qed
\enddemo

From Theorem~2.5 and the above result we obtain

\proclaim{Corollary~4.3}
Let \mvp U\ be a univalent, bounded degree \cp\ with $\carr U =\bold C$.
Denote by \mcompk\ the complex of \mvp U .
Then the set of all finite valence \cp s for \mcompk\ with finite branch sets is equal to the set of all range packings of \dcp s for \mcompk.
Moreover, each such a packing is determined uniquely, up to similarities of $\bold C$, by its branch set.
\endproclaim

\remark{Concluding Remarks}
(1)~We conjecture that Theorem~4.2 should also be true for triangulations with recurrent 1-skeletons but not necessarily of bounded degree.

(2)~For complete analogy with the univalent setting, one should be able to show the uniqueness result, similar to Corollary~4.3, for branched \cp s whose corresponding univalent packings cover $\bold C$ but underlying complexes are not of bounded degree.
As it was shown in [HSc], such complexes might have 1-skeletons which are not recurrent even though univalent packings associated with them have carriers equal to $\bold C$.
\endremark

\enddocument